\documentclass[11pt,twoside]{amsart}
\usepackage{amssymb}
\newtheorem{theorem}{Theorem}[section]

\newtheorem{lemma}[theorem]{Lemma}
\newtheorem{corollary}[theorem]{Corollary}

\theoremstyle{definition}
\newtheorem{defn}[theorem]{Definition}

\newtheorem{remark}[theorem]{Remark}

\newcommand{\skipit}[1]{{}}
\newcommand{\prfend}{\hbox to7pt{\hfil}
\par\vskip-\baselineskip\hbox to\hsize
{\hfil\vbox {\hrule width6pt height6pt}}\vskip\baselineskip}

\newcommand{\N}{\mathbb{N}}

\newcommand{\myarrow}[2]{\hbox to #1pt{\hfil$\to$\hfil}{\hskip-#1pt{\raise
10pt\hbox to#1pt{\hfil$\scriptscriptstyle #2$\hfil}}}}

\begin{document}

\title{Simplicial complexes and Macaulay's inverse systems}

\author{Adam Van Tuyl}
\address{Department of Mathematical Sciences, Lakehead University, Thunder Bay, ON P7B 5E1, Canada}
\email{avantuyl@lakeheadu.ca}

\author{Fabrizio Zanello}
\address{Department of Mathematical Sciences, Michigan Technological University, Houghton, MI 49931-1295, USA}
\email{zanello@math.kth.se}


\begin{abstract}
Let $\Delta$ be a simplicial complex on $V = \{x_1,\ldots,x_n\}$, with Stanley-Reisner ideal  
$I_{\Delta}\subseteq R = k[x_1,\ldots,x_n]$. The goal of this paper is to investigate the class of 
artinian algebras $A=A(\Delta,a_1,\ldots,a_n)= R/(I_{\Delta},x_1^{a_1},\ldots,x_n^{a_n})$, where each 
$a_i \geq 2$.  By utilizing the technique of Macaulay's inverse systems, we can explicitly describe the
socle of $A$ in terms of $\Delta$. As a consequence, we determine the simplicial complexes, that we will 
call {\em levelable}, for which there exists a tuple $(a_1,\ldots,a_n)$ such that $A(\Delta,a_1,\ldots,a_n)$ 
is a level algebra.
\end{abstract}

\maketitle

\section{Introduction}

The purpose of this paper is to study a class of monomial artinian algebras naturally associated with any given simplicial complex. Our approach will introduce the use of {\em Macaulay's inverse systems}, a tool coming from commutative algebra, to this subject. In particular, we characterize which simplicial complexes are {\em levelable}, that is, when at least one of the associated artinian algebras is level.

We now summarize the results of this paper (and postpone formal definitions until later sections).
Let $\Delta$ be a simplicial complex on the vertex set $V = \{x_1,\ldots,x_n\}$ with associated Stanley-Reisner 
ring $R/I_{\Delta}$.  By adjoining the pure powers of the variables to $I_{\Delta}$, we can make a monomial 
artinian algebra from $R/I_{\Delta}$.  Specifically, if $(a_1,\ldots,a_n) \in \N^n$ with $a_i \geq 2$ for all $i$, set
\[A(\Delta,a_1,\ldots,a_n) := R/(I_{\Delta}+(x_1^{a_1},\ldots,x_n^{a_n})).\]

Our paper will extend past work on artinian algebras of the form $A(\Delta,a_1,\ldots,a_n)$, which have  mainly been studied in the special case that $(a_1,a_2,\ldots,a_n) = (2,2,\ldots,2)$. For example, Aramova, Herzog and 
Hibi \cite{AHH}, who called these algebras {\it indicator algebras}, gave a formula for the Betti numbers of 
$A$ when $I_{\Delta}$ is a squarefree strongly stable ideal. Herzog and Hibi \cite{HH} used information about 
the socle of $A$ to give an upper-bound on the number of facets of the simplicial complex $\Delta$; and when 
$I_{\Delta} = I(G)$ is the edge ideal of a graph, Simis, Vasconcelos, and Villarreal \cite{SVV} related the  
type of $A$ to the number of minimal vertex covers of $G$, under certain hypotheses on $G$. In fact, we will 
prove a similar result (see Corollary \ref{applicationgraph}). Hibi's paper \cite{H} is one of the first papers to study the connections between level algebras and simplicial complexes, and has since inspired a significant amount of research by both Hibi and other authors (see, e.g., \cite{HH2},\cite{HHV},\cite{H2}).

The first half of this note is devoted to describing the {\em socle}, i.e., the annihilator of the maximal 
homogeneous ideal, of the artinian algebra $A=A(\Delta,a_1,\ldots,a_n)$.  We show (see Theorem \ref{geninverse}) that 
the socle-vector of $A$ is a function of the {\em facets} (the maximal faces) of $\Delta$ and the  tuple 
$(a_1,\ldots,a_n)$. Our proof of this fact is based upon  {\em Macaulay's inverse systems}. 
In fact, we give an explicit description of the generators of the inverse system. Moreover, we also prove the 
interesting fact that the number of minimal generators of the inverse system, or equivalently, the Cohen-Macaulay type, of $A$ is independent of the $a_i$'s, and is exactly the number of facets of $\Delta$. The second part of the paper then focuses on the problem of determining for which $a_i$'s, if any,  $A$ is a {level} algebra.

Level algebras, which were introduced in commutative algebra by Stanley in his seminal paper \cite{St}, have been (and continue to be) the focus of an active area of research.  The memoir \cite{GHMS} is perhaps the best source of references for this subject up to the year 2003. For more on level algebras, especially for the applications outside commutative algebra or combinatorics, such as to invariant theory, geometry, and even complexity theory (the connection being the theory of inverse systems),  see the references \cite{IK}, \cite{KS}, and \cite{NW}, among others.  

Using our description of the socle,  we prove (see Theorem \ref{maintheorem}) that $A(\Delta,a_1,\ldots,a_n)$ is 
level if and only if $(a_1,\ldots,a_n)$ is a solution to a system of linear equations of a special type, where 
these equations are constructed from the facets of $\Delta$.  Our result enables one to utilize computer algebra 
programs, such as Mathematica, to determine if one can construct a level algebra starting from $\Delta$.

We then introduce the notion of a {\it levelable simplicial complex} to denote a simplicial complex 
$\Delta$ for which $A(\Delta,a_1,\ldots,a_n)$ is level for some tuple $(a_1,\ldots,a_n)$.  After proving that 
nonlevelable simplicial complexes exist, we show that some interesting families of simplicial complexes, 
such as simplicial forests (see Theorem  \ref{algore}), are levelable.

It is hoped that this note will encourage future work on levelable simplicial complexes, since our paper has just begun investigating this subject.  Many interesting questions still deserve to be addressed.

\section{Preliminaries}

In this section we collect together the needed definitions and results about simplicial complexes and level 
artinian algebras for the latter part of the paper. We consider standard graded algebras $A=R/I$, 
where $R=k[x_1,...,x_n]$, $I$ is a homogeneous ideal of $R$, and $\deg x_i = 1$ for all $i$.   For 
simplicity, we will assume that the field $k$ has characteristic zero (although many of the results 
of this paper can be proved for any infinite field). 

The algebra $A$ is called {\em artinian} if its Krull-dimension is 0. We then have that $A$ is artinian 
if and only if $\sqrt{I}$, the radical of $I$, equals $(x_1,\ldots,x_n)$, the maximal homogeneous ideal of $R$,
if and only if the Hilbert function of $A$ is eventually 0. Hence, in the artinian case, we can simply speak of the 
{\em $h$-vector} of $A$, $h(A)=h=(h_0,h_1,...,h_e)$, where $h_i=\dim_k A_i$ and $e$ is the last index such 
that $\dim_k A_e>0$. Since we may suppose without loss of generality that $I$ does not contain non-zero forms 
of degree 1, $h_1=n$ (the number of variables of $R$) is defined as the {\em codimension} of $A$.

Let $\Delta$ be a {\em simplicial complex} on the vertex set $V = \{x_1,\ldots,x_n\}$.  That is, $\Delta$
is a collection of subsets of $V$ such that: $(i)$ if $F \in \Delta$ and $G \subset F$, then $G \in \Delta$; $(ii)$ for all $i = 1,\ldots,n$, $\{x_i\} \in \Delta$.  We let $\mathcal{F}(\Delta ) = \{F_1,\ldots,F_t\}$ denote the elements
of $\Delta$ that are maximal under inclusion;  these elements are called the {\em facets} of $\Delta$,
while any element of $\Delta$ is called a {\em face}. The complex $\Delta$ is  {\em pure} if all its facets have the same cardinality. 

We associate to $\Delta$ the squarefree monomial ideal $I_{\Delta}$ in the polynomial ring $R = k[x_1,\ldots,x_n]$ defined by $ I_{\Delta} := (\{ x_{i_1}\cdots x_{i_s} ~|~ \{x_{i_1},\ldots,x_{i_s}\} \not\in \Delta\}).$ The ideal $I_{\Delta}$, which is generated by the set of non-faces of $\Delta$, is called the {\em Stanley-Reisner ideal}
of $\Delta$, and $R/I_{\Delta}$ is the {\em Stanley-Reisner ring}. 
Properties of $\Delta$ are then encoded into the algebraic invariants of $R/I_{\Delta}$;  see \cite{BH,SB} for more. 

Except when $\Delta = \emptyset$,  the ring $R/I_{\Delta}$ is never an artinian algebra, since $\sqrt{I_{\Delta}} = I_{\Delta}\neq (x_1,...,x_n)$. However, by adjoining powers of the variables, we can construct an artinian algebra from $R/I_{\Delta}$.  In particular,
\[A = A(\Delta,a_1,\ldots,a_n) = R/(I_{\Delta},x_1^{a_1},\ldots,x_n^{a_n}). ~~\]
Without loss of generality, we can assume that $a_i \geq 2$ for all $i$. Indeed, if some $a_i = 1$, then we can reduce to 
\[A \cong k[x_1,\ldots,\hat{x}_i,\ldots,x_n]/(I_{\Delta'},x_1^{a_1}, \ldots,\hat{x}_i,\ldots,x_n^{a_n}) =A(\Delta',a_1,\ldots,\hat{a}_i,\ldots,a_n),\]
where $\Delta'$ is the simplicial complex with face set $\{W \subseteq \Delta ~|~ W \subseteq V\setminus\{x_i\}\}$.
A similar argument also allows us to assume that $\Delta$ contains no facet of the form $\{x_i\}$.
We therefore make the convention for the remainder of the paper that $a_i \geq 2$ 
for all $i$ and that $\Delta$ has no facets of cardinality one.

\section{The socle of $A(\Delta,a_1,\ldots,a_n)$}

The goal of this section is to describe the socle-vector of the artinian algebras $A(\Delta,a_1,\ldots,a_n)$ by using
Macaulay's inverse system.  The {\em socle} of $A$ is the annihilator of the maximal homogeneous ideal 
$(\overline{x}_1,...,\overline{x}_n)$ of $A$, that is, soc$(A)=\lbrace a\in A {\ } \mid {\ } a\overline{x}_i=\overline{0}$ for all $i\rbrace $. Since soc$(A)$ is a homogeneous ideal, we can define the {\em socle-vector} of $A$ as 
$s(A)=s=(s_0,s_1,...,s_e)$, where $s_i=\dim_k$ soc$(A)_i$. Notice 
that $s_0=0$ and $s_e=h_e>0$. The integer $e$ is called the {\em socle degree} of $A$ (or of $h$). 
The {\em (Cohen-Macaulay) type} of the socle-vector $s$ (or of the algebra $A$) is type$(s)=\sum_{i=0}^es_i$. In particular, if $s=(0,0,...,0,s_e=t)$, we say that the algebra $A$ (or its $h$-vector) is {\em level} (of type $t$). 
Moreover, if $t=1$, then $A$ is said to be {\em Gorenstein}.

Let us now introduce {\em {(Macaulay's)} inverse systems}, also known as {\em Matlis duality},
which will play a key r\^ole in this note. For a complete introduction to this theory, please see \cite{Ge} and \cite{IK}. 

Given the polynomial ring $R=k[x_1,x_2,...,x_n]$, let $S:=k[y_1,y_2,...,y_n]$, and consider $S$ as a graded $R$-module, 
where the action of $x_i$ on $S$ is partial differentiation with respect to $y_i$. There is a one-to-one correspondence between artinian algebras $A=R/I$ and finitely generated $R$-submodules $M$ of $S$, where $I=\operatorname{Ann}(M)$ is the annihilator of $M$ in $R$ and, conversely, $M=I^{-1}$ is the $R$-submodule of $S$ which is annihilated by $I$ (cf. \cite{Ge}, Remark 1, p.17).

If $R/I$ has socle-vector $s=(s_0,s_1,...,s_e)$, then $M$ is minimally generated by $s_i$ elements of degree $i$, 
for $i=1,2,...,e$, and the $h$-vector of $R/I$ is given by the number of linearly independent partial derivatives 
obtained in each degree by differentiating the generators of $M$ (cf. \cite{Ge}, Remark 2, p.17). In particular, 
level algebras of type $t$ and socle degree $e$ correspond to $R$-submodules of $S$ minimally generated by $t$ 
elements of degree $e$.

It can be shown that the inverse system  $M$ of an ideal $I$ is generated by monomials if and only if $I$ is a 
monomial ideal (and this is the case in which we are interested in this note). Precisely, if $I$ is monomial, the 
monomials in $M$ in each degree ``are'' exactly those not in $I$ in the same degree (after, of course, 
renaming the $x_i$'s with the corresponding $y_i$'s).   

Our first main result shows that the minimal generators of the inverse system corresponding to 
$A(\Delta,a_1,\ldots,a_n)$ are directly obtained from the facets of $\Delta$, and that, therefore, 
the socle-vector of such an algebra can be explicitly determined.  We begin with the basic case $A(\Delta,2,\ldots,2)$:

\begin{lemma}\label{is}
If $\Delta$ is a simplicial complex with facet set $\mathcal{F}(\Delta)$, then the minimal generators of the 
inverse system corresponding to $A(\Delta,2,\ldots,2)$ are the elements of the set
\[\left.\left\{\prod_{x_i \in F} y_i ~\right|~ F \in \mathcal{F}(\Delta)\right\}.\]
\end{lemma}

\begin{proof}
We have seen above that the inverse system $M$ of a monomial ideal $I$ contains exactly the monomials  not 
in $I$ rewritten in the variables $y_i$'s. Moreover, since we operate on the module $M$ by partial  
differentiation, it is easy to see that a minimal set of generators for $M$ is given by those monomials of $M$ having no non-constant multiple in $M$. Therefore, because the  ideal we are considering here is $I=(I_{\Delta},x_1^2,...,x_n^2)$, it immediately follows that the facets of ${\Delta}$ give the desired generators, as stated.
\end{proof}

\begin{theorem}\label{geninverse}
Let $\Delta$ be a simplicial complex having facet set $\mathcal{F}(\Delta)$.  Let us fix the algebra $A =A(\Delta,a_1,\ldots,a_n)$, with $a_i \geq 2$ for all $i$.  Then the minimal generators of the inverse system corresponding to $A$ are the elements of the set
\[\left.\left\{\prod_{x_i \in F} y_i^{a_i-1} ~\right|~ F \in \mathcal{F}(\Delta)\right\}.\]
In particular, if $s(A)= (s_0,\ldots,s_e)$ is the socle-vector of $A$, then
\[s_j = \#\left\{F \in \mathcal{F}(\Delta) ~\left|~ \sum_{x_i \in F} (a_i-1) = j\right\}\right.. \]
\end{theorem}

\begin{proof} Again, we want to find all the ``maximal'' monomials not in 
$I:=(I_{\Delta},x_i^{a_1},\ldots,x_n^{a_n})$ (expressed in the $y_i$'s). We first prove that the monomials 
$\prod_{x_i \in F} y_i^{a_i-1}$ belong to $M$ and  are maximal there, for every $F \in \mathcal{F}(\Delta)$. 
Indeed, similarly to what we have observed in the proof of Lemma \ref{is}, those monomials must be in $M$ since 
they correspond to facets; the fact that they are maximal is because their exponents are the highest possible, since by adding one to any of them we would end up in $I$.

It remains to show that those are the only maximal monomials in $M$. But since, by inverse systems, 
any squarefree monomial dividing a monomial of $M$ does not belong to $I$ and therefore corresponds to 
a face of $\Delta $, it is immediate to see that any other generator of $M$ would necessarily imply the 
existence of a new facet of $\Delta $, which is a contradiction. 

The last statement is an immediate consequence of the relationship between generators of the inverse system and entries of the socle-vector.
\end{proof}

\begin{remark}
Herzog and Hibi \cite{HH} proved this result in the special case that
$A = A(\Delta,2\ldots,2)$ by describing the canonical module $\omega_A$ of $A$.
The inverse system approach highlights the combinatorial nature of the socle-vector of $A$.
\end{remark}

\begin{corollary} \label{Gorenstein}
Let $\Delta$ be a simplicial complex, and let $A = A(\Delta,a_1,\ldots,a_n)$, 
for any choice of the $a_i \geq 2$.  Then:
\begin{enumerate}
\item[i)] The Cohen-Macaulay type of  
$A$ equals  $|\mathcal{F}(\Delta)|$, and is therefore independent of the $a_i$'s;
\item[ii)] $A(\Delta,a_1,\ldots,a_n)$  is Gorenstein if and only $\Delta $ has a unique facet.
\end{enumerate}
\end{corollary}

As an application of  Theorem \ref{geninverse}, the socle-vector can be used to obtain information about
independent sets of {\em finite simple graphs} $G = (V_G,E_G)$ (i.e., there are no loops or multiple
edges) on the vertex set $V_G = \{x_1,\ldots,x_n\}$ and with edge set $E_G$.

Recall that the {\em edge ideal} of $G$ is the ideal $I(G):= (\{ x_ix_j ~|~ \{x_i,x_j\} \in E_G\} ) \subseteq R = k[x_1,\ldots,x_n].$  We then say that $W \subseteq V_G$ is an {\em independent set} if no two vertices of $W$ are adjacent in $G$.

\begin{corollary} \label{applicationgraph}
Let $G$ be a graph on the
vertex set $V = \{x_1,\ldots,x_n\}$, and let $I(G) \subseteq  R = k[x_1,\ldots,x_n]$ be the corresponding edge ideal. If
\[\mathbb{F}: 0 \rightarrow  F_{n}=R^{d_1}(-c_1) \oplus \cdots \oplus R^{d_t}(-c_t) \rightarrow F_{n-1} \rightarrow \cdots \rightarrow F_1 \rightarrow R \rightarrow A  \rightarrow 0 \]
is the graded minimal free resolution of $A= R/(I(G),x_1^2,\ldots,x_n^2)$, 
then $d_i$ is the number of maximal independent sets of $G$ of size  $c_i - n$.  In particular,
\[\operatorname{type}(A) = d_1 + \cdots + d_t = \mbox{number of maximal independent sets.}\]
\end{corollary}

\begin{proof}
For any artinian algebra $A = R/I$ with socle-vector $s(A) = (s_0,\ldots,s_e)$ and codimension $n$, a {\em graded minimal free resolution} (MFR) of $A$ is a long exact sequence of  the form
\[\mathbb{F}: 0 \rightarrow F_n=\bigoplus_{1 \leq i \leq e} R^{s_i}(-i-n) \rightarrow F_{n-1} \rightarrow \cdots \rightarrow F_1 \rightarrow R \rightarrow A \rightarrow 0. \]
It is then well known that the {\em graded Betti numbers} (i.e., the multiplicities of the shifts) of the last module of the MFR of $A$ are given by the entries of the socle-vector of $A$. 
In particular, $s_i = 0$ if and only if  $R^{s_i}(-i-n)$ does not appear.

Since $I(G)$ is a squarefree monomial ideal, it is also the Stanley-Reisner ideal of some simplicial complex $\Delta$. Specifically, $I(G)$ corresponds to the simplicial complex  
\[\Delta = \Delta_G := \{W =\{x_{i_1},\ldots,x_{i_s}\} \subseteq V_G ~|~ \mbox{$W$ is an independent set of $G$}\}.\]
In other words, $I(G)$ is the Stanley-Reisner ideal of the simplicial complex $\Delta_G$ whose facets correspond to maximal independent sets of $G$.  
Because of the relationship  between the MFR and the socle-vector, the MFR of $R/(I(G),x_1^2,\ldots,x_n^2)$
encodes information about the independent sets, and thus the first conclusion now follows from Theorem \ref{geninverse}.  

The second statement is an immediate consequence of the definition of type.
\end{proof}

\section{Levelable simplicial complexes}

We now focus on studying which algebras $A(\Delta,a_1,\ldots,a_n)$, associated to the simplicial complex 
$\Delta $, are level. Thanks to Theorem \ref{geninverse}, we can already prove the main result of this section. 
Note that in the theorem below we exclude the case that $t=1$, because by Corollary \ref{Gorenstein},
$A(\Delta,a_1,\ldots,a_n)$ is level (in fact, Gorenstein) for all choices of the $a_i$'s in this case.

\begin{theorem}\label{maintheorem}
Let $\Delta$ be a simplicial  complex on $n$ vertices with facets $\mathcal{F}(\Delta) = \{F_1,\ldots,F_t\}$, where $t \geq 2$. Let $F_i = \{x_{i,1},\ldots,x_{i,{d_i}}\}$ denote the $i$th facet.  Then the algebra $A(\Delta,a_1,\ldots,a_n)$, with each $a_i\geq 2$, is level if and only if $(a_1,\ldots,a_n)$ is a simultaneous integral solution to the following $t-1$ equations:
\begin{eqnarray*}
(x_{1,1} + \cdots + x_{1,d_1}) - (x_{2,1}+\cdots + x_{2,d_2})  & = & d_1-d_2\\
(x_{2,1} + \cdots + x_{2,d_2}) - (x_{3,1}+\cdots + x_{3,d_3})  & = & d_2-d_3 \\
& \vdots & \\
(x_{t-1,1} + \cdots + x_{t-1,d_{t-1}}) - (x_{t,1}+ \cdots + x_{t,d_t}) & = & d_{t-1} -d_{t}.
\end{eqnarray*}
\end{theorem}

\begin{proof} Denote a facet of $\Delta$ by $F_i = \{x_{i,1},\ldots,x_{i,d_i}\}$. The algebra 
$A(\Delta,a_1,\ldots,a_n)$, with socle-vector $s=(s_0,...,s_e)$, is level if and only if $s_j = 0$ 
for all $j \neq e$, if and only if, by Theorem \ref{geninverse},
\[\sum_{k=1}^{d_1} (a_{1,k}-1) = \sum_{k=1}^{d_2} (a_{2,k}-1) = \cdots = \sum_{k=1}^{d_t} (a_{1,k}-1).\]
But this is equivalent to:
\begin{eqnarray*}
a_{1,1} + \cdots + a_{1,d_1} - d_1 &=& a_{2,1} + \cdots + a_{2,d_2} - d_2\\
a_{2,1} + \cdots + a_{2,d_2} - d_2&=&  a_{3,1} + \cdots + a_{3,d_3} - d_3\\
& \vdots &\\ 
a_{t-1,1} + \cdots + a_{t-1,d_{t-1}} - d_{t-1}&=&  a_{t,1} + \cdots + a_{t,d_t} - d_t,
\end{eqnarray*}
and the theorem follows.
\end{proof}

\begin{remark}\label{alex}
Determining whether or not a solution to the above system exists for any given simplicial complex can be accomplished
on many computer algebra systems.  Michigan Tech University grad students Matt Miller and Alex Schaefer 
(the latter being a student of the second author) have written one such program in Mathematica, which can 
be found on the first author's webpage.\footnote[1]{http://flash.lakeheadu.ca/$\sim $avantuyl/research/Levelable\_VanTuyl\_Zanello.html}
\end{remark}

Theorem \ref{maintheorem} immediately implies a result found in Boij's Thesis \cite{B}:

\begin{corollary} \label{mats}  Let $\Delta$ be a simplicial complex.
\begin{enumerate}
\item[i)] {\em (Boij)} $A(\Delta,2,\ldots,2)$ is level if and only if $\Delta$ is pure;
\item[ii)] More generally, if $\Delta$ is {pure}, then $A(\Delta,d,\ldots,d)$ is level for 
every  $d \geq 2$.
\end{enumerate}
\end{corollary}

If there exists a solution to the equations of Theorem \ref{maintheorem}, we give $\Delta$ the following name:

\begin{defn}\label{levelable}
A simplicial complex $\Delta $ is {\em levelable} if there exists a tuple $(a_1,\ldots,a_n)$, with each $a_i\geq 2$, such that the artinian algebra $A(\Delta,a_1,\ldots,a_n)$ is level.
\end{defn}

We now show that the set of tuples making $\Delta $ levelable, if non-empty, is infinite:

\begin{corollary}\label{infinite}
If $\Delta $ is levelable for some choice of $(a_1,\ldots,a_n) \in \N^n$, with 
$a_i \geq 2$, then it is levelable for infinitely many such choices.
\end{corollary}

\begin{proof}
Suppose there exists an $(a_1,\ldots,a_n) \in \N^n$, with $a_i \geq 2$, such that $A(\Delta,a_1,\ldots,a_n)$ is level. 
Then a straightforward calculation shows that, for every positive integer $c$, $(c(a_1-1)+1,\ldots,c(a_n-1)+1)$ 
also satisfies the system of equations of Theorem \ref{maintheorem}. Therefore all the 
algebras $A(\Delta,c(a_1-1)+1,\ldots,c(a_n-1)+1)$ are also level.
\end{proof}

Not all simplicial complexes are levelable. Indeed, we have:

\begin{theorem}\label{nonlev} For every $n \geq 5$, there exists a  simplicial complex on $n$ vertices which is not levelable.
\end{theorem}

\begin{proof}
We first consider the case that $n\geq 5$ and $n$ is odd.  Let $\Delta_n$ denote the simplicial complex whose facet set is given by
\[ \mathcal{F}(\Delta_n) = \{\{x_1,x_3,x_5,\ldots,x_n\},\{x_2,x_4,x_6,\ldots,x_{n-1}\}, \{x_1,x_4,x_6,\ldots,x_{n-1}\}, \{x_2,x_5,x_7,\ldots,x_n\}\}.\]
The first facet has $(n+1)/2$ elements, while the other facets have $(n-1)/2$ elements.

Now consider the system
\begin{eqnarray*}
(x_1+x_3+\cdots+x_n) - (x_2+x_4+\cdots+x_{n-1}) &=& (n+1)/2 - (n-1)/2 = 1 \\
(x_2+x_4 +\cdots+x_{n-1}) - (x_1+x_4+\cdots+x_{n-1}) &=& (n-1)/2 - (n-1)/2 = 0 \\
(x_1+x_4+ \cdots+x_{n-1}) - (x_2+x_5+\cdots+x_n) & = & (n-1)/2 - (n-1)/2 = 0.
\end{eqnarray*}

The second equation implies that $x_1 = x_2$.  Combining this result with the 
third equation gives $(x_4+\cdots+x_{n-1}) = (x_5+\cdots+x_n)$.  By subbing these identities into the first equation 
we obtain $x_3 =1$.  Thus, any integral solution $(a_1,\ldots,a_n)$ to the above system must have $a_3 =1$, and therefore, 
by Theorem \ref{maintheorem}, $\Delta_n$ is not levelable.

When $n\geq 5$ and $n$ is even, we consider the simplicial complex $\Delta_n$ with facet set
\[\mathcal{F}(\Delta_n) = \{\{x_1,x_3,x_5,x_6,x_7,\ldots,x_n\},\{x_2,x_5,x_6,x_7,\ldots,x_n\},\{x_1,x_4\}, \{x_2,x_4\}\}.\]
In a similar way, one can show (we leave the routine details to the reader) that $\Delta_n$ is not
levelable, because any integral solution to the system defined by Theorem \ref{maintheorem} would again have $a_3 = 1$.
\end{proof}

We round out this paper by showing that some interesting classes of simplicial complexes are 
levelable.  In particular, we show that the assumption of $n\geq 5$ was necessary in constructing the
non-levelable simplicial complexes in Theorem \ref{nonlev}. 

\begin{theorem} The following simplicial complexes are levelable:
\begin{enumerate}
\item[i)] Any pure simplicial complex;
\item[ii)] Any simplicial complex on $n \leq 4$ vertices;
\item[iii)] Any simplicial complex $\Delta$ with pairwise disjoint facets.
\end{enumerate}
\end{theorem}

\begin{proof} i) This fact follows immediately follows from Corollary \ref{mats}.

For statement ii), one can simply check all non-isomorphic simplicial complexes on
four or less vertices.  (Note that i) already implies that many of these simplicial complexes are levelable.)

In order to show  iii), suppose that $\mathcal{F}(\Delta) = \{F_1,\ldots,F_t\}$.  The proof is by induction on $t$. If $t=1$, then $\Delta$ is clearly levelable by Corollary \ref{Gorenstein}.   When $t=2$, suppose $F_1 = \{x_1,\ldots,x_{d_1}\}$ and $F_2 = \{y_1,\ldots,y_{d_2}\}$.  Setting $x_2 = \cdots = x_{d_1}=y_2 = \cdots = y_{d_2} = 2$,
we consider the equation
\[x_1 +  2(d_1-1) - y_1 - 2(d_2-1) = d_1 -d_2 \Leftrightarrow x_1 = d_1-d_2 + 2(d_2-d_1) + y_1.\]
We can now find an integer $a \geq 2$ such that $d_1-d_2 + 2(d_2-d_1) + a \geq 2$.
Hence $\Delta$ will be levelable by Theorem \ref{maintheorem} if we take
the tuple $(d_2-d_1+a,2,\ldots,2,a,2,\ldots,2)$. 

Suppose now that $\mathcal{F}(\Delta) = \{F_1,\ldots,F_t\}$ where  $t > 2$.  By
induction, the simplicial complex $\Delta'$ with facets $\mathcal{F}(\Delta')= \{F_1,\ldots,F_{t-1}\}$ is levelable.  We can assume that $\Delta'$ is a simplicial complex on $\{x_1,\ldots,x_m\}$, 
and $(a_1,\ldots,a_m)$ is a tuple that makes $\Delta'$ levelable.

If $\Delta$ is a simplicial complex on the vertex set $\{x_1,\ldots,x_m,x_{m+1},\ldots,x_n\}$,
then because the facets are disjoint, $F_t = \{x_{m+1},\ldots,x_n\}$.  After relabeling,
we can assume that $F_{t-1} = \{x_{p+1},\ldots,x_m\}$.  Now consider the equation
\[(x_{p+1} + \cdots + x_m) - (x_{m+1} + \cdots + x_n) = (m-p) - (n-m) = 2m -n - p.\]
Set $x_{m+1} = \cdots = x_{n-1} = 2$, and rearrange to give
\[x_n = -n+p+2 +(x_{p+1}+\cdots+ x_m).\]
We can now pick an integer $c$ sufficiently large such that
\[e:= -n+p+2+ (c(a_{p+1}-1)+1 + \cdots + c(a_m-1)+1) \geq 2.\]
If $c$ is such an integer, then $\Delta$ is levelable since the tuple $(c(a_1-1)+1,\ldots,c(a_m-1)+1,2,\ldots,2,e)$ satisfies the system of Theorem \ref{maintheorem}.  Indeed, we have just shown that this tuple will satisfy the last equation of the system.  The remaining equations will only involve $x_1,\ldots,x_m$.  As shown in the
proof of Theorem \ref{infinite}, the tuple $(c(a_1-1)+1,\ldots,c(a_m-1)+1)$
will satisfy the first $t-1$ equations of Theorem \ref{maintheorem}.
\end{proof}

Let us now turn our attention to another class of complexes.  Faridi \cite{F} introduced simplicial forests, an important subclass of simplicial complexes that generalizes the notion of a forest from
graph theory.  To define a simplicial forest, one must begin by defining a leaf.

\begin{defn}  Let $F$ be a facet of a simplicial complex $\Delta$. Then $F$ is a {\em leaf} if either $F$ is the only facet of $\Delta$, or there exists some facet $G \neq F$ such that for all facets $G' \in \Delta$ with $G' \neq F$, we have $G' \cap F \subseteq G \cap F.$
We then call $\Delta$ a {\em forest} if every simplicial complex generated by a sub-collection of the facets of $\Delta$ has a leaf.  When $\Delta$ is connected, we say $\Delta$ is a {\em tree}.
\end{defn}

It is well known (see, for example, \cite{F}, Remark 2.3) that if $F$ is a leaf of $\Delta$, then $F$ contains at least one vertex that does not belong to any other facet of $\Delta$.  We can now prove:

\begin{theorem}\label{algore}  If $\Delta$ is a simplicial forest, then $\Delta$ is levelable.\footnote[2]{As a humorous aside, perhaps a good name for this theorem is the ``Clearcutting Theorem'', since we have unwittingly
proved that all forests are levelable!}
\end{theorem}

\begin{proof} The proof is by induction on the number of facets of $\Delta$.  
When $|\mathcal{F}(\Delta)|=1$, the result follows  since $A$ is Gorenstein.  
So,  suppose $\mathcal{F}(\Delta) = \{F_1,F_2\}$. By relabeling the vertices, we can assume
\[F_1 = \{x_1,\ldots,x_p,x_{p+1},\ldots,x_m\} ~~ \mbox{and} ~~ F_2 = \{x_{p+1},\ldots,x_m,x_{m+1},\ldots,x_n\}.\]
According to Theorem \ref{maintheorem}, we need to find  an integral solution $(a_1,\ldots,a_n)$ to 
\[(x_1+\cdots+x_p + x_{p+1} + \cdots +x_m)-(x_{p+1}+\cdots+x_m+x_{m+1}+\cdots+x_n) = m-(n-p),\]
such that $a_i \geq 2$.  Setting $x_2 = \cdots = x_{n-1} = 2$ and rearranging gives
\[x_1 = n-m - p + x_n.\]
If $e: = n-m - p  \geq 0$, the desired solution is $(e+2,2,\ldots,2)$.  If $e < 0$, then we can pick any $a_n \geq 2$ such that $e+a_n \geq 2$.  The desired solution is then $(e+a_n,2,\ldots,2,a_n)$.   

Suppose now that $\Delta$ has $t$ facets.  Since $\Delta$ is a forest,
let $F \in \Delta$ be a leaf of $\Delta$, and let $G$ be the 
facet of $\Delta$ such that $F \cap G' \subseteq F \cap G$ for every
facet $F \neq G' \in \Delta$.  If $\mathcal{F}(\Delta) = \{F_1,\ldots,F_t\}$
are the facets of $\Delta$, we relabel so that $F_{t-1} = G$ and $F_t = F$. Furthermore, after relabeling the vertices, we can assume that
\[F_{t-1} = \{x_p,\ldots,x_s,x_{s+1},\ldots,x_m\} ~~\mbox{and}~~ F_t = \{x_{s+1},\ldots,x_m,x_{m+1},\ldots,x_n\}.\]
Note that the vertices $\{x_{m+1},\ldots,x_n\}$ only appear in the facet
$F_t$, while the remaining vertices of $F_t$ also appear in $F_{t-1}$. 

If $\Delta'$ is the simplicial complex generated by the facets
$\{F_1,\ldots,F_{t-1}\}$, then $\Delta'$ is also a simplicial forest
on the vertex set $\{x_1,\ldots,x_m\}$, and thus, by induction, there exists a tuple $(a'_1,\ldots,a'_m)$ such
that $A(\Delta',a'_1,\ldots,a'_m)$ is level. We now consider the equation
\begin{equation}\label{oneeq}
(x_p+\cdots + x_s + x_{s+1} + \cdots + x_m) -  (x_{s+1}+\cdots+x_m+x_{m+1}+\cdots + x_n) 
\end{equation}
\[= (x_p + \cdots +x_s) - (x_{m+1} + \cdots + x_n) = (m-(p-1))-(n-s).\]
Setting $x_{m+1}=\cdots = x_{n-1} = 2$ and rearranging gives
\[ (x_p + \cdots +x_s) -n+m-s+p+1= x_n.\]
Now pick any integer $c \geq 1$ such that
\[ e:=(c(a'_p-1)+1) + \cdots + (c(a'_s-1)+1)-n+m-s+p+1   \geq 2.\]
Clearly such a $c$ exists.  
Then, we claim that the desired tuple is
\[(c(a'_1-1)+1,\ldots,c(a'_m-1)+1,2,\ldots,2,e).\]
Indeed, we have just shown that this tuple satisfies (\ref{oneeq}). On the other hand, we know by Corollary \ref{infinite} that $(c(a'_1-1)+1,\ldots,c(a'_m-1)+1)$ also satisfies all the equations of Theorem \ref{maintheorem} for $\Delta'$.  But the equations associated to $\Delta$ are precisely the equations of $\Delta'$, plus equation $(\ref{oneeq})$.  Since $x_{m+1},\ldots,x_n$ do not appear in any of these equations, the tuple $(c(a'_1-1)+1,\ldots,c(a'_m-1)+1,2,\ldots,2,e)$ is  therefore a solution, for $(c(a'_1-1)+1,\ldots,c(a'_m-1)+1)$ is a solution for all the other equations.
\end{proof}
\noindent
{\ }\\
{\bf Acknowledgements.} The work of this paper began when the first author visited the second author at Michigan Tech University.  The first author would like to thank Michigan Tech for its hospitality.  He would also like to acknowledge the financial support of NSERC while working on this project. 

We warmly thank our former supervisor, Tony Geramita, and Tony Iarrobino for helpful comments. We also want to thank Michigan Tech grad students Matt Miller and Alex Schaefer for writing for us a computer program in Mathematica (see Remark \ref{alex}). Many of the computations and examples produced for this paper were performed with the help of CoCoA \cite{Co}.


\end{document}